\documentclass{amsart}
\usepackage{amsmath}
\usepackage{amssymb}
\usepackage{amscd}
\usepackage{xcolor}
\usepackage{hyperref}
\hypersetup{hidelinks}

\newtheorem{theorem}{Theorem}[section]
\newtheorem{thm}[theorem]{Theorem}
\newtheorem{cor}[theorem]{Corollary}
\newtheorem{lem}[theorem]{Lemma}
\newtheorem{prop}[theorem]{Proposition}

\theoremstyle{definition}
\newtheorem{defn}[theorem]{Definition}

\newtheorem{rem}[theorem]{Remark}
\newtheorem{constr}[theorem]{Construction}

\newtheorem{conj}[theorem]{Conjecture}

\theoremstyle{remark}

\newcommand{\mbb}{\mathbb}
\newcommand{\QQ}{\mbb{Q}}

\newcommand{\ZZ}{\mbb{Z}}
\newcommand{\CC}{\mbb{C}}

\newcommand{\AAA}{\mbb{A}}
\newcommand{\PP}{\mbb{P}}

\newcommand{\mc}{\mathcal}

\newcommand{\mcT}{\mc{T}}

\newcommand{\mcX}{\mc{X}}

\newcommand{\mfm}{\mathfrak{m}}

\newcommand{\OO}{\mc{O}}

\newcommand{\wht}{\widehat}

\newcommand{\whtOO}{\wht{\OO}}

\newcommand{\SP}{\text{Spec }}

\newsavebox{\sembox}
\newlength{\semwidth}
\newlength{\boxwidth}

\newsavebox{\semrbox}
\newlength{\semrwidth}
\newlength{\boxrwidth}

\newcommand{\obhat}{\whtOO_{B, b}}

\newcommand{\xbhat}{\wht{\mcX}_b}

\title
[Integral Hodge/Tate conjecture of quotients via weak approximation]{Integral Hodge/Tate conjecture of quotient varieties via weak approximation}

\author[Cao]{Yang Cao}
\address{
School of Mathematics\\
Shandong University\\
250100, Jinan, China}
\email{yangcao1988@gmail.com}

\author[Tian]{Zhiyu Tian}
\address{
Beijing International Center for Mathematical Research\\
Peking University\\
100871, Beijing, China}
\email{zhiyutian@bicmr.pku.edu.cn}

\date{\today}

\begin{document}


\begingroup\color{black}
\begin{abstract}
For a finite group acting on a smooth projective separably rationally connected variety over an algebraically closed field, we prove that the rational quotient map induces a surjection from the coinvariants of its second Kato homology onto the second Kato homology of any smooth projective quotient model (with coefficients invertible in the ground field). The integral Hodge and Tate conjectures for one-cycles consequently descend to such quotient models. The proof uses weak approximation over function fields to retain the image of a smooth fiber in the total space.
\end{abstract}
\endgroup

\maketitle

\section{Introduction}



The main motivation of this paper is the following conjecture proposed by Voisin.
\begin{conj}[Voisin]\label{conj:IHC_RC}
 Let $X$ be a smooth projective separably rationally connected variety of dimension $d$ defined over an algebraically closed field $k$. For every prime $\ell\ne\operatorname{char}k$, the group $H^{2d-2}_{\text{\'et}}(X,\ZZ_\ell(d-1))$ is algebraic. If $X$ is defined over the complex numbers, then $H^{2d-2}_{\text{sing}}(X,\ZZ)$ is algebraic.
\end{conj}
Here algebraicity with $\ZZ_\ell$ coefficients means surjectivity of the cycle class map after tensoring $CH_1$ with $\ZZ_\ell$. The integral assertion over $\CC$ includes torsion classes. Throughout the paper, $\ell$-adic and finite \'etale coefficients are understood to be invertible in the ground field.

This conjecture is commonly referred to as \emph{the integral Hodge or Tate conjecture for one-cycles} on (separably) rationally connected varieties.

We remark that the standard formulation of the Tate or Hodge conjecture asserts that the cycle class map surjects onto a subspace of the cohomology group which is either invariant under a suitable Galois group (in the Tate case) or belongs to the $(p,p)$-component of the Hodge decomposition (in the Hodge case).
For (separably) rationally connected varieties, thanks to the existence of a decomposition of the diagonal \`a la Bloch--Srinivas \cite{BlochSrinivas}, it is known that every class in $H^{2d-2}_{\text{\'et}}(X,\QQ_\ell(d-1))$ or $H^{2d-2}_{\text{sing}}(X,\QQ)$ is algebraic.
Consequently, the integral version asks whether the entire integral cohomology space is spanned by classes of algebraic cycles.

The most compelling evidence supporting Conjecture~\ref{conj:IHC_RC} is as follows. Let us first recall the Tate conjecture for varieties defined over a finite field.

\begin{conj}[Tate conjecture]\label{conj:Tate}
 Let $X$ be a smooth projective variety of dimension $d$ defined over a finite field. Then the cycle class map
 \[
 CH^i(X)\otimes\QQ_\ell\longrightarrow H^{2i}_{\text{\'et}}(X,\QQ_\ell(i))
 \]
 is surjective.
\end{conj}
For $X/\mathbb F_q$, the target in Conjecture~\ref{conj:Tate} is arithmetic \'etale cohomology. By Hochschild--Serre and weights it identifies with
\[
H^{2i}_{\text{\'et}}(X_{\overline{\mathbb F}_q},\QQ_\ell(i))^{\operatorname{Gal}(\overline{\mathbb F}_q/\mathbb F_q)}.
\]

Schoen \cite{SchoenIntegralTate} established that the classical formulation of the Tate conjecture for divisor classes on surfaces (with $\QQ_\ell$-coefficients) leads to an integral version of the Tate conjecture for curve classes on all smooth projective varieties over an algebraic closure of a finite field.
By employing deformation theory of rational curves together with Schoen's result, Voisin obtained the following:

\begin{thm}[{\cite[Theorem 1.6]{VoisinCurveClassRC}}]\label{thm:Voisin}
 If the Tate conjecture (Conjecture~\ref{conj:Tate}) holds for divisor classes on surfaces, then Conjecture~\ref{conj:IHC_RC} holds for smooth projective rationally connected varieties over $\CC$.
\end{thm}

Some special cases were previously known.

\begin{thm}\label{thm:IHC}
 Conjecture~\ref{conj:IHC_RC} is true for rationally connected smooth projective varieties defined over an algebraically closed field of characteristic $0$ in the following cases:
 \begin{enumerate}
 \item (Voisin \cite{VoisinIHC3fold}) $X$ has dimension $3$.
 \item (H\"oring--Voisin \cite{HoeringVoisin}) $X$ is a smooth Fano $4$-fold.
 \end{enumerate}
\end{thm}

\begingroup\color{black}
When the first draft of this paper was written, Conjecture~\ref{conj:IHC_RC} seemed out of reach.
However, with the recent development of AI, the landscape has changed. 
The second named author has noticed that Voisin's conjecture \ref{conj:IHC_RC} is now a theorem, whose proof is due to AI with some human assistance.

\begingroup\color{black}
In this paper, mostly motivated by Voisin's conjecture \ref{conj:IHC_RC}, we study how the obstruction to algebraicity of curve classes changes under finite group quotients. Our main result concerns the second Kato homology, in the homological conventions of \cite[Example 2.2 and Definition 2.3]{zerocycleLaurent}. For a smooth projective variety this group is the cokernel of the cycle class map. With the Borel--Moore homology convention recalled in Section~\ref{sec:kato}, the finite-coefficient sequence is
\[
 CH_1(X)\otimes_\ZZ\Lambda
 \longrightarrow H_2(X,\Lambda)\longrightarrow KH_2(X,\Lambda)
 \longrightarrow0,\qquad\Lambda=\ZZ/n.
\]
 Over $\CC$, the integral singular-homology sequence is
\[
 CH_1(X)\longrightarrow H_2(X,\ZZ)
 \longrightarrow KH_2(X,\ZZ)\longrightarrow0.
\]

\begin{thm}\label{thm:main}
Let $k$ be an algebraically closed field. Let $X/k$ be smooth projective and separably rationally connected, and
let a finite group $G$ act on $X$. Suppose $Y$ is a smooth projective
birational model of $X/G$. For every positive integer $n$ invertible in
$k$, the rational quotient map $q:X\dashrightarrow Y$ induces a
surjection
\begin{equation}\label{eq:main}
 KH_2(X,\ZZ/n)_G\twoheadrightarrow KH_2(Y,\ZZ/n).
\end{equation}
Here the subscript denotes coinvariants for the action by proper
pushforward. 
If $k=\CC$, there is
also a surjection
\begin{equation}\label{eq:integral-main}
 KH_2(X,\ZZ)_G\twoheadrightarrow KH_2(Y,\ZZ)
\end{equation}
for the singular Borel--Moore homology convention of
\cite[Example 2.2(1)]{zerocycleLaurent}.
\end{thm}

The map in \eqref{eq:main} is the pushforward defined by the graph correspondence and is compatible with smooth birational changes of the quotient model, as proved in Lemma~\ref{lem:graphs}. The subscript $G$ means the quotient by the subgroup generated by $g_*\xi-\xi$, for $g\in G$. The hypothesis that $Y$ is a smooth projective model is part of the statement. If $|G|$ is invertible in $k$, tame resolution supplies such a model; see Remark~\ref{rem:G-resolution}.

\endgroup

Theorem~\ref{thm:main} implies that the validity of Voisin's conjecture \ref{conj:IHC_RC} descends to finite group quotients, and in particular establishes it for smooth projective models of tame finite quotients of projective space.

\begin{cor}\label{cor:quotientITC}
Let $X$ be a smooth projective separably rationally connected variety over an algebraically closed field $k$, let a finite group $G$ act on $X$, and let $Y$ be a smooth projective birational model of $X/G$. If Conjecture~\ref{conj:IHC_RC} holds for $X$, then it holds for $Y$. In particular, it holds for every smooth projective model of $\PP^d/G$ when $|G|$ is invertible in $k$; such models exist in this tame case.
\end{cor}

\begin{rem}\label{rem:G-resolution}
Even though resolution of singularities is unknown in positive characteristic, a variety with tame quotient singularities over a perfect field admits a projective resolution, functorial for smooth morphisms; see \cite[Theorem E]{BerghRydhFunctorialDestackification}. In particular, if $H\triangleleft G$ acts tamely on a smooth projective variety $V$, there is a smooth projective resolution
\begin{equation}\label{eq:equivariant-resolution}
 \mathcal R(V/H)\longrightarrow V/H
\end{equation}
carrying the residual $G/H$-action. This is the equivariance needed in Proposition~\ref{prop:finite}. More generally, every smooth projective model $V'$ of a finite quotient of a separably rationally connected $V$ is separably rationally connected: the induced dominant rational map $V\dashrightarrow V'$ has a separable function-field extension, and a general member of a separably dominating family of rational curves avoids its codimension-at-least-two indeterminacy locus. Composing this family with the rational map retains separable dominance of its two-point evaluation. This last assertion needs a supplied smooth model, but not a tame action.
\end{rem}

An interesting feature of the proof of Theorem~\ref{thm:main} is that we use an arithmetic property of rationally connected varieties over function fields: weak approximation. We refer readers to Section~\ref{sec:WA} for its formulation. The following theorem is the link between weak approximation and Kato homology.

\begin{thm}\label{thm:WAKH2}
Let $\pi:\mcX\to B$ be a flat projective morphism of relative dimension $d\geq1$, with $\mcX$ smooth projective and $B$ a smooth projective curve over an algebraically closed field $k$. Suppose the geometric generic fiber is smooth separably rationally connected. Assume that for every finite surjective morphism $C\to B$ from a smooth projective curve, weak approximation holds for $\mcX\times_B C\to C$. Let $n$ be invertible in $k$ and put $\Lambda=\ZZ/n$. For every smooth fiber $i_c:\mcX_c\hookrightarrow\mcX$,
\begin{equation}\label{eq:fibre-image}
\begin{split}
 H^{2d}(\mcX,\Lambda(d))={}&\operatorname{Im}\bigl(CH_1(\mcX)/n\xrightarrow{\operatorname{cl}}H^{2d}(\mcX,\Lambda(d))\bigr)\\
 &+i_{c,*}H^{2d-2}(\mcX_c,\Lambda(d-1)).
\end{split}
\end{equation}
Consequently the proper pushforward induces a surjection
\begin{equation}\label{eq:fibre-kato}
 i_{c,*}:KH_2(\mcX_c,\Lambda)\twoheadrightarrow KH_2(\mcX,\Lambda).
\end{equation}
\end{thm}
\endgroup

\begingroup\color{black}
Theorem~\ref{thm:WAKH2} is proved in Section~\ref{sec:fibration} using local one-cycles and the Mayer--Vietoris sequence. In Section~\ref{sec:quotients}, a cyclic isotrivial family identifies \eqref{eq:fibre-kato} with the actual quotient map by a cyclic group. Induction and Sylow subgroup descent then prove Theorem~\ref{thm:main}. The consequence of Theorem~\ref{thm:WAKH2} for the integral Hodge and Tate conjectures is Corollary~\ref{cor:WAITC}.
\endgroup

\textbf{Acknowledgement and AI usage:} The proof strategy and most parts of the proof were obtained by the authors near the end of the year 2025, and a sketchy draft was finished around March 2026 but left incomplete until recently. In the first draft, only the implication: IHC/ITC of $X$ $\implies$ IHC/ITC for the quotient was stated. But Z.T. has another argument using arithmetic properties over a Laurent field to prove the comparison of the relevant Kato homologies. Recently, AI was used to check the draft and Z.T. asked AI to check if the argument can be strengthened to prove the stronger comparison result. AI found some inaccuracies of the original draft and fixed them, also verified that the stronger comparison follows from the same line of argument.

Z.T. would like to thank the AI for Math team at PKU, in particular, Liang Xiao and Shuai Chen, for providing the necessary AI tools.

Both Z.T. and Y.C. are partially supported by the National Key Project for Research and Development Grant No. 2025YFA1017300. Z.T. is partially supported by NSFC Grant No. 12525102.

\section{Preliminaries}
\subsection{Weak approximation}\label{sec:WA}
In this section, we recall the formulation of weak approximation.
Let $k$ be an algebraically closed field, and let $B$ be a smooth projective curve over $k$, with function field $K=k(B)$. We assume $X_K$ is a smooth projective variety over $K$.

Given a sequence $(\nu_1,\ldots,\nu_l)$ of finitely many places of $K$, each $\nu_i$ corresponds to a closed point $b_i\in B$. Let $\whtOO_{B,b_i}$ denote the completion of the local ring $\OO_{B,b_i}$ at the maximal ideal $\mfm_{B,b_i}$, and let $\wht K_{b_i}$ be its fraction field. Then $X_K$ satisfies weak approximation over $(\nu_1,\ldots,\nu_l)$ if the image of the diagonal map
\[
X_K(K)\longrightarrow\prod_i X_K(\wht K_{b_i})
\]
is dense, where the target takes the product of the $\nu_i$-adic topologies.

Now we will unwind the definition of weak approximation using a geometric formulation. For a detailed discussion, one may find \cite[Section 1]{WASurvey} very helpful.

\begin{defn}\label{def:model}
A \emph{model} of $X_K$ is a flat projective morphism $\pi:\mcX\to B$ with generic fiber isomorphic to $X_K$.
\end{defn}
\noindent Each section of the model corresponds to a $K$-rational point of $X_K$, and vice versa, by properness and the regularity of $B$.

Now let us fix a model $\pi:\mcX\to B$, a place $\nu_i$, and a section $\wht s_i$ of the restriction
\[
\pi|_{\wht B_{b_i}}:\mcX\times_B\wht B_{b_i}\longrightarrow\wht B_{b_i},
\qquad \wht B_{b_i}=\SP\whtOO_{B,b_i}.
\]
\noindent Note that $\wht s_i$ corresponds to a point of $X_K(\wht K_{b_i})$.

For each nonnegative integer $N$, a basic neighborhood of $\wht s_i$ consists of sections of this restriction congruent to $\wht s_i$ modulo $\mfm_{B,b_i}^{N+1}$. Then one can reformulate weak approximation in the following way.
\begin{defn}\label{defn:wa}
Let $X_K$ be a smooth projective variety over $K=k(B)$, where $k$ is algebraically closed. For a fixed model $\pi:\mcX\to B$, say that approximation holds \emph{at order $N$ on this model} if, for every finite sequence of distinct closed points $b_1,\ldots,b_l$ and formal sections $\wht s_i$ over $\wht B_{b_i}$, there is a section $\sigma:B\to\mcX$ agreeing with $\wht s_i$ modulo $\mfm_{B,b_i}^{N+1}$ for each $i$.
We say that $X_K$ satisfies \emph{weak approximation} if this holds at every order $N\geq0$ on one model.
Once a model $\mcX\to B$ is fixed, we also say that this family satisfies weak approximation if its generic fiber does.
\end{defn}
Full weak approximation only depends on the generic fiber, so the condition at all orders in Definition~\ref{defn:wa} holds for any model if it holds for one model. A change of model can change the numerical orders defining an individual neighborhood; no model-independence of one fixed order $N$ is asserted.
\begin{rem}
    At the time when a first draft was prepared, weak approximation seemed out of reach. 
    However, with the recent development of AI, the situation has changed. 
The second named author has been informed that at least over the complex numbers, weak approximation can be proven using a recent AI generated preprint, which is complex analytic in nature.
In general, in a joint work in preparation, the second named author and Jason Starr give a proof of the full weak approximation for any separably rationally connected variety as defined in \ref{defn:wa}.
\end{rem}

\subsection{\texorpdfstring{\textcolor{black}{Kato homology and rational maps}}{Kato homology and rational maps}}\label{sec:kato}

\begingroup\color{black}
Let $k$ be algebraically closed, let $n>0$ be invertible in $k$, and put $\Lambda=\ZZ/n$. We use Borel--Moore homology $H_a$ and Kato homology $KH_a$ as in \cite[Example 2.2 and Definition 2.3]{zerocycleLaurent}. For a smooth projective $r$-fold $V$ and a finite coefficient module $A$ of order invertible in $k$, the underlying homology is
\[
 H_a(V,A)=H^{2r-a}_{\mathrm{et}}(V,A(r)).
\]
The low-degree exact sequence \cite[equation (3.1)]{zerocycleLaurent}, with Tate twists displayed, is
\begin{equation}\label{eq:defect}
 (CH_1(V)/n)\otimes_\Lambda\Lambda(1)
 \xrightarrow{\operatorname{cl}}H_2(V,\Lambda)
 \to KH_2(V,\Lambda)\to0.
\end{equation}
Equivalently, the usual cycle class map has target
\[
 H^{2r-2}_{\mathrm{et}}(V,\Lambda(r-1))=H_2(V,\Lambda(-1))
\]
and gives
\begin{equation}\label{eq:kato-twist}
\begin{split}
 KH_2(V,\Lambda)&\simeq
 \operatorname{coker}\bigl(CH_1(V)/n\xrightarrow{\operatorname{cl}}
 H_2(V,\Lambda(-1))\bigr)\otimes_\Lambda\Lambda(1).
\end{split}
\end{equation}
Here $\Lambda(1)$ is a constant invertible rank-one $\Lambda$-module, so tensoring by it is exact and commutes with proper pushforward and coinvariants. Quotienting cycles by algebraic equivalence does not alter their image. Over $\CC$ we also use singular Borel--Moore homology and its associated Kato homology with integral coefficients, as in \cite[Example 2.2(1)]{zerocycleLaurent}. The groups $H_2$ and $KH_2$ vanish for zero-dimensional varieties.
\endgroup

\begingroup\color{black}
\begin{lem}\label{lem:boundary}
Let $S,T$ be smooth projective varieties of dimensions $s,t$, and
let $R\in CH_s(S\times T)$ be represented by a cycle supported on
$D\times T$ for a proper closed subset $D\subsetneq S$.
Then $R_*$ maps $H_2(S,\Lambda)$ into the image of the cycle class map
\[
 CH_1(T)\otimes_\ZZ \Lambda\longrightarrow H_2(T,\Lambda).
\]
\end{lem}
\begin{proof}
The Chinese remainder theorem reduces the proof to the case $n=\ell^a$. 
Choose an irreducible component of $R$ with reduced structure
and take $D$ to be its first projection image, of dimension $e<s$.
A smooth projective alteration $D'\to D$ of degree $b$ prime to
$\ell$ exists by \cite[Expos\'e X, Th\'eor\`eme 2.1]{IllusieLaszloOrgogozoTravauxGabber};
this is also the proper-support argument in
\cite[proof of Lemma 2.36]{Hasse0cycle}.
Pull the component back over the finite flat locus and take the cycle
closure $R'\in CH_s(D'\times T)$. With $j:D'\to S$,
\begin{equation}\label{eq:boundary-action}
 (j\times1)_*R'=bR,\qquad bR_*\alpha=R'_*j^*\alpha.
\end{equation}
If $e<s-1$, the group $H^{2s-2}(D',\Lambda(s))$ vanishes.
If $e=s-1$, it is top cohomology, generated by the class of a rational point  tensored with
$\Lambda$; the algebraic correspondence $R'$ carries these classes
into the image of the cycle class map on $T$. Invertibility of $b$ and
\eqref{eq:boundary-action} prove the assertion. 
\end{proof}

\begin{lem}\label{lem:graphs}
A dominant rational map $f:S\dashrightarrow T$ between smooth
projective varieties induces a pushforward
$f_*:KH_2(S,\Lambda)\to KH_2(T,\Lambda)$. These maps are functorial for dominant
rational maps. Birational maps and the projection
$T\times\PP^1\to T$ induce isomorphisms. For a generically finite
map of degree $m$, its transpose graph induces $f^*$ with
\begin{equation}\label{eq:pushpull}
 f_*f^*=m\operatorname{id}.
\end{equation}
Moreover, if $i:F\hookrightarrow S$ is a closed immersion of smooth
projective varieties, $f$ is defined generically along $i(F)$, and
$q=f|_F:F\dashrightarrow T$ is dominant, then
\begin{equation}\label{eq:restriction}
 f_*i_*=q_*:KH_2(F,\Lambda)\longrightarrow KH_2(T,\Lambda).
\end{equation}
\end{lem}
\begin{proof}
If $\Gamma_f\subset S\times T$ is the graph closure, its action is
\begin{equation}\label{eq:graph-action}
 (\Gamma_f)_*\alpha=(p_T)_*
 \bigl(p_S^*\alpha\cup\operatorname{cl}(\Gamma_f)\bigr).
\end{equation}
The dimensions and twists in \eqref{eq:graph-action} give an element
of $H_2(T,\Lambda)$, and the image of the cycle class map is preserved.
For dominant $f,g$, the correspondence
$\Gamma_g\circ\Gamma_f-\Gamma_{g\circ f}$ is supported over a proper
closed subset of $S$: on the dense open where both maps are morphisms
it is zero, with graph multiplicity one. Chow localization and
Lemma~\ref{lem:boundary} give functoriality, hence birational invariance.
The same argument with source $F$ proves \eqref{eq:restriction}, since
$f$ is a morphism near the generic point of $i(F)$.

The projective bundle formula decomposes $H_2(T\times\PP^1,\Lambda)$ into
the summand from $T\times\{a\}$ and a top-cohomology summand generated
by the classes of $\{x\}\times\PP^1$ tensored with $\Lambda(1)$.
The latter is algebraic and is killed by
projection. Projection is the identity on the first summand, giving
the asserted isomorphism.
Hence the projection $T\times\PP^1\to T$ induces an isomorphism.

For a generically finite $f$, finite flat push--pull on a dense open
$V\subset T$ and Chow localization give
\[
 \Gamma_f\circ{}^t\Gamma_f=m\Delta_T+R,
 \qquad\operatorname{Supp}R\subset(T\setminus V)\times T.
\]
The transpose preserves curve classes because the dimensions agree.
Lemma~\ref{lem:boundary} proves \eqref{eq:pushpull}.

If a smooth graph resolution $u:S'\to S$, $v:S'\to T$ is available,
functoriality gives $f_*=v_*u_*^{-1}$. Thus this definition agrees
with resolution-defined pushforward and is independent of choices.
\end{proof}
\endgroup

\begin{lem}\label{lem:reduction_mod_l}
Let $X$ be a smooth projective variety of dimension $d$ over an algebraically closed field $k$, and let $\ell\ne\operatorname{char}k$. Surjectivity of
\[
 CH_1(X)/\ell\longrightarrow H^{2d-2}(X,\ZZ/\ell(d-1))
\]
implies surjectivity of
\[
 CH_1(X)\otimes\ZZ_\ell\longrightarrow H^{2d-2}(X,\ZZ_\ell(d-1)).
\]
If $H^{2d-1}(X,\ZZ_\ell(d-1))$ is torsion free, the converse also holds.
\end{lem}
\begin{proof}
The coefficient sequence always supplies an injection
\begin{equation}\label{eq:coefficient-injection}
 H^{2d-2}(X,\ZZ_\ell(d-1))/\ell
 \hookrightarrow H^{2d-2}(X,\ZZ/\ell(d-1)).
\end{equation}
Every finite cycle class has an integral cycle lift. If finite cycle surjectivity holds, \eqref{eq:coefficient-injection} shows that the cokernel of the $\ell$-adic cycle map is zero modulo $\ell$. This cokernel is a finitely generated $\ZZ_\ell$-module, so it is zero by Nakayama's lemma. Under the stated torsion-freeness assumption, \eqref{eq:coefficient-injection} is an isomorphism, proving the converse.
\end{proof}

\section{\texorpdfstring{\textcolor{black}{Weak approximation and second Kato homology}}{Weak approximation and second Kato homology}}\label{sec:fibration}

Throughout this section, we fix the following notations.
Fix an algebraically closed field $k$.
Let $n$ be a positive integer, invertible over the field $k$, and put $\Lambda=\ZZ/n$.
All cohomology groups are \'etale cohomology unless stated otherwise explicitly; the Tate twists are displayed.
Let $\pi:\mcX\to B$ be a flat projective morphism of relative dimension $d\geq1$ from a smooth projective variety to a smooth projective curve, with smooth geometrically connected generic fiber.

Denote by $U\subset B$ a nonempty open subset over which $\pi$ is smooth. We may shrink $U$ when necessary.
For a point $b\in B$, let $R_b=\obhat$ be the completion of $\OO_{B,b}$ and let $K_b$ be its fraction field.
Denote by $\xbhat$ the base change of $\mcX$ to $\SP R_b$, and put $\xbhat^\circ=\mcX\times_B\SP K_b$.
Given any morphism $T\to B$, denote by $\mcX_T$ the base change to $T$.
The completed models $\xbhat$ are regular, because completion of the local ring of the smooth base curve is a regular morphism.

\subsection{Analysis of one-cycles}

\begin{lem}\label{lem:trait-generators}
Let $R=k[[t]]$ and let $\mathcal Y$ be a regular projective flat $R$-scheme of pure relative dimension $d\geq1$ with smooth generic fiber. 
Write $\mathcal Y_0=\sum_\lambda m_\lambda D_\lambda$ the special fibre. For a closed smooth point $x\in D_\lambda$ on no other component, there is a closed subscheme $T_\lambda\subset\mathcal Y$ which is a regular finite flat $R$-curve, generically \'etale over $R$, and satisfies
\begin{equation}\label{eq:trait-intersections}
 \deg(T_\lambda/R)=m_\lambda,\qquad
 (T_\lambda\cdot D_\mu)=\delta_{\lambda\mu}.
\end{equation}
\end{lem}
\begin{proof}
Fix a closed immersion $\mathcal Y\hookrightarrow\PP_R^M$.
Choose $d$ hyperplanes through $x$ on $\PP_k^M$ whose common intersection with the special fiber has finite support and whose differentials restrict to a basis of $T_x^*D_\lambda$.

Lift their linear equations to $\ell_1^{(0)},\ldots,\ell_d^{(0)}$ over $R$. For a fixed $N\geq2$, allow all coefficient perturbations
\begin{equation}\label{eq:trait-lifts}
 \ell_j=\ell_j^{(0)}+t^Nq_j,\qquad
 q_j\in H^0(\PP_R^M,\OO(1)),\quad 1\leq j\leq d.
\end{equation}
These preserve the reductions and first jets at $x$.
Put $K=k((t))$. The coefficient tuples in \eqref{eq:trait-lifts} form a coset $a+t^NR^{d(M+1)}$, which is Zariski dense in $\AAA_K^{d(M+1)}$: a polynomial over $K$ vanishing on $R$ is zero because $R$ is infinite; induction proves density of $R^s$, and translation and multiplication by $t^N$ preserve it.
Generic Bertini on the smooth embedded $K$-scheme $\mathcal Y_K$ gives a nonempty open set of tuples with smooth zero-dimensional intersection \cite[Lemma 33.47.3]{StacksProjectBertiniTheorems}. Thus a tuple in \eqref{eq:trait-lifts} can be chosen in this open. The scheme
\[
 Z=\mathcal Y\cap V(\ell_1,\ldots,\ell_d)
\]
has finite \'etale generic fiber and finite closed fiber. It is finite over $R$ by projectivity \cite[Lemma 30.21.2]{StacksProjectApplicationsFormalFunctions}.

Set $A=\OO_{\mathcal Y,x}$ and choose an equation $u$ of $D_\lambda$. As $x$ lies on no other component,
\begin{equation}\label{eq:trait-parameters}
 \begin{gathered}
 t=v u^{m_\lambda},\qquad v\in A^\times,\\
 u,\ell_1,\ldots,\ell_d\text{ is a regular system of parameters in }A.
 \end{gathered}
\end{equation}
Here the local equations of the hyperplanes are taken in a trivialization of $\OO(1)$. Their restrictions give a cotangent basis in $A/(u)=\OO_{D_\lambda,x}$, proving the second assertion of \eqref{eq:trait-parameters}. Hence $S=A/(\ell_1,\ldots,\ell_d)=\OO_{Z,x}$ is a DVR with uniformizer $\bar u$.
Since $R$ is complete henselian, $\SP S$ is the open-and-closed local factor of the finite $R$-scheme $Z$ at $x$ \cite[Lemma 10.153.4(1)]{StacksProjectHenselianLocalRings}. Take this factor as $T_\lambda$; it is closed in $\mathcal Y$.

By \eqref{eq:trait-parameters}, $t$ is nonzero in $S$, so $S$ is torsion-free and hence flat over $R$. Its generic fiber is a factor of $Z_K$, hence \'etale. Its residue field is $k$, and
\[
 \operatorname{rank}_R S=\dim_k(S/tS)
 =\operatorname{length}_S(S/(\bar u^{m_\lambda}))=m_\lambda.
\]
The valuation of $u$ is one, and no other special-fiber component meets this local factor. This proves \eqref{eq:trait-intersections}. Separability is imposed by Bertini on the generic fiber, not inferred from the ramification index.
\end{proof}

\begin{lem}\label{lem:local-classes}
Write $\mcX_b=\sum_\lambda m_\lambda\Lambda_\lambda$ as a fiber divisor. There is an isomorphism
\begin{equation}\label{eq:component-trace}
 \theta_b:H^{2d}(\xbhat,\Lambda(d))\xrightarrow{\sim}
 \bigoplus_\lambda\Lambda,\qquad
 \beta\longmapsto\left(\operatorname{Tr}_{\Lambda_\lambda}
 (\beta|_{\Lambda_\lambda})\right)_\lambda.
\end{equation}
Every class in this group is the class of $h_*[T]$ for a finite flat normal generically \'etale $R_b$-scheme $T$ and an $R_b$-map $h:T\to\xbhat$.
\end{lem}
\begin{proof}
Proper base change identifies the group with the top cohomology of the reduced special fiber. Restriction to the irreducible components is an isomorphism in degree $2d$, because their intersections have dimension at most $d-1$. The top cohomology of each integral projective component is $\Lambda$, with its degree trace. One may see the latter assertion by normalization and removal of a subset of dimension at most $d-1$. This proves \eqref{eq:component-trace}; compare \cite[Appendix A, (A.5)]{WittenbergEsnault_0_cycle}. This cohomological assertion requires neither normal crossings nor characteristic zero.

Apply Lemma~\ref{lem:trait-generators} to the regular projective $R_b$-model $\xbhat$, with $R_b\simeq k[[t_b]]$. For each $\Lambda_\lambda$, choose a closed smooth point on no other component. Equation~\eqref{eq:trait-intersections} supplies a regular finite flat generically \'etale branch $T_\lambda$ of degree $m_\lambda$ and intersection vector $\delta_{\lambda\mu}$.
Compatibility of Cartier intersection, cycle classes and proper pushforward gives
\begin{equation}\label{eq:local-intersections}
 \theta_b\bigl(\operatorname{cl}(h_*[T])\bigr)_\lambda
 =\deg h^*\Lambda_\lambda\pmod n.
\end{equation}
For a prescribed class, represent its coordinates in \eqref{eq:component-trace} by $a_\lambda\in\{0,\ldots,n-1\}$ and take the disjoint union of $a_\lambda$ copies of these branches. Equations~\eqref{eq:trait-intersections} and \eqref{eq:local-intersections} prove the assertion, with degree $\sum_\lambda a_\lambda m_\lambda$. The source and its proper pushforward retain all multiplicities.
\end{proof}

\begin{lem}\label{lem:map-continuity}
Let $T$ be a finite flat normal generically \'etale $R_b$-scheme and $h:T\to\xbhat$ an $R_b$-map. There is an integer $N$ such that every $R_b$-map $g:T\to\xbhat$ agreeing with $h$ modulo $t_b^N$ satisfies
\[
 \operatorname{cl}(g_*[T])=\operatorname{cl}(h_*[T])
 \quad\text{in }H^{2d}(\xbhat,\Lambda(d)).
\]
\end{lem}
\begin{proof}
On each complete DVR branch $S_a$ of $T$, put $e_a=\operatorname{ord}_{S_a}(t_b)$. For each component through its closed image choose a local equation $f_\lambda$, and set $r_{a\lambda}=\operatorname{ord}_{S_a}(h^*f_\lambda)$. These orders are finite because $h$ is over the trait. Take $N$ with $Ne_a>r_{a\lambda}$ for every pair. Congruence modulo $t_b^N$ preserves the closed images, and both maps land in the chosen local neighborhoods: an open subset of a local scheme containing its closed point is the whole scheme. The pullbacks of $f_\lambda$ then differ by an element of $t_b^NS_a$, so have the same order. Components not through the closed image have unit pullbacks. The result follows from \eqref{eq:local-intersections} and the injectivity in \eqref{eq:component-trace}.
\end{proof}

\begin{lem}\label{lem:globalization}
Let $B$ be a smooth projective curve and $b_1,\ldots,b_m\in B$ distinct closed points. Let $Z_i$ be finite flat normal schemes of the same positive degree $D$ over $R_i=\whtOO_{B,b_i}$, generically \'etale over $R_i$. Then there is a smooth integral projective curve $C$ and a finite flat generically \'etale morphism $q:C\to B$ of degree $D$ such that
\begin{equation}\label{eq:globalized-covers}
 C\times_B\SP R_i\simeq Z_i\quad\text{over }R_i,\qquad 1\leq i\leq m.
\end{equation}
\end{lem}
\begin{proof}
Write $Z_i=\SP A_i$, with $A_i=\prod_a S_{ia}$ a product of complete DVRs. If $s_{ia}$ is a uniformizer and $e_{ia}$ is the rank over $R_i$, then $t_i$ is a unit times $s_{ia}^{e_{ia}}$. The elements $1,s_{ia},\ldots,s_{ia}^{e_{ia}-1}$ give a basis modulo $t_i$, hence an $R_i$-basis by Nakayama. Thus $S_{ia}=R_i[s_{ia}]$, with monic relation reducing to $T^{e_{ia}}$. Shift the generators by constants of $k$ distinct for the different factors. Their monic relations then have pairwise coprime reductions, hence unit resultants, so the Chinese remainder theorem gives
\[
 A_i=R_i[T]/P_i(T),\qquad P_i\text{ monic of degree }D.
\]

We apply Hensel's lemma as follows. Write $u_a$ for the component of the image of $T$ in a branch $S_a$ of $A_i$, and put $e_a=v_a(t_i)$ and $\delta_a=v_a(P_i'(u_a))<\infty$. Choose $N$ such that $Ne_a>2\delta_a$ and $Ne_a-\delta_a\geq e_a$ for every branch. If a monic degree-$D$ polynomial $Q$ satisfies $Q\equiv P_i\pmod{t_i^N}$, then
\[
 v_a(Q(u_a))\geq Ne_a,\qquad v_a(Q'(u_a))=\delta_a.
\]
The strong Hensel lemma produces roots $u'_a\in S_a$ with $v_a(u'_a-u_a)\geq Ne_a-\delta_a\geq e_a$. Thus the induced map
\begin{equation}\label{eq:algebra-stability}
 R_i[T]/Q\longrightarrow A_i,\qquad T\longmapsto(u'_a)_a
\end{equation}
is an isomorphism modulo $t_i$, hence surjective by Nakayama, and is an isomorphism because both modules are free of rank $D$. This is the application of Hensel's lemma needed in \eqref{eq:globalized-covers}; compare the smooth-branch statement of \cite[Lemma 1.4(b)]{deJongStarr_GHS}.

Choose an auxiliary point $q_0\notin\{b_1,\ldots,b_m\}$ and an affine open $V\subset B$ containing all these points. By the Chinese remainder theorem, approximate the coefficients of the $P_i$ to the orders required by the application of Hensel's lemma in \eqref{eq:algebra-stability}, while imposing
\[
 P(T)\equiv T^D+\tau\pmod{\mfm_{q_0}^2},
 \qquad \operatorname{ord}_{q_0}(\tau)=1.
\]
The resulting monic polynomial $P\in\Gamma(V,\OO_V)[T]$ is Eisenstein at $q_0$, hence irreducible over $k(B)$. It is separable because its completed algebra at each $b_i$ is generically \'etale. Let $C$ be the normalization of $B$ in $k(B)[T]/P$. It is a normal integral projective curve, hence smooth over the perfect field $k$. Its finite morphism to $B$ is flat of degree $D$, because it is torsion-free over the local DVRs. Normalization commutes with completion for these excellent curves. Since each prescribed $A_i$ is already normal, \eqref{eq:algebra-stability} identifies these completed normalizations with $A_i$, proving \eqref{eq:globalized-covers}. If there are no marked points, add one prescribed \'etale algebra of rank $D$ and apply the preceding argument. For $D=1$, one may simply take $C=B$.
\end{proof}

\subsection{Analysis of cohomology}
\begin{prop}\label{prop:MVsequence}
Assume $B\setminus U$ is finite. There is a long exact sequence containing
\begin{equation}\label{eq:MVsequence}
\begin{split}
 \bigoplus_{b\notin U}H^{2d-1}(\xbhat^\circ,\Lambda(d))
 &\xrightarrow{\partial}H^{2d}(\mcX,\Lambda(d))\\
 &\longrightarrow H^{2d}(\mcX_U,\Lambda(d))\oplus
 \bigoplus_{b\notin U}H^{2d}(\xbhat,\Lambda(d))\\
 &\longrightarrow\bigoplus_{b\notin U}H^{2d}(\xbhat^\circ,\Lambda(d)).
\end{split}
\end{equation}
The same comparison applies to a constructible $\Lambda$-complex on the base.
\end{prop}
\begin{proof}
For a closed point $b$ and its complement $V=B\setminus\{b\}$, localization gives the comparison
\begingroup\small
\[
\begin{CD}
 H^{j-1}(\mcX_V,\Lambda(d)) @>>> H^j_{\mcX_b}(\mcX,\Lambda(d)) @>>> H^j(\mcX,\Lambda(d)) @>>> H^j(\mcX_V,\Lambda(d))\\
 @VVV @VV\simeq V @VVV @VVV\\
 H^{j-1}(\xbhat^\circ,\Lambda(d)) @>>> H^j_{\mcX_b}(\xbhat,\Lambda(d)) @>>> H^j(\xbhat,\Lambda(d)) @>>> H^j(\xbhat^\circ,\Lambda(d)).
\end{CD}
\]
\endgroup
To justify the support isomorphism, apply localization on $B$ to $R\pi_*\Lambda(d)$ and excise to the henselian trait at $b$. Passing from that excellent henselian trait to its completion preserves the closed-stalk cohomology by proper base change and the punctured-trait cohomology by the equivalence of finite separable extensions of their fraction fields. The localization triangles therefore identify the support groups. Comparing these triangles gives the sequence for one point; excising all points gives \eqref{eq:MVsequence}. The same reasoning applies to constructible complexes on $B$.
\end{proof}

\begin{lem}\label{lem:surj}
Denote $f_U:\mcX_U\to U$. Let $F\in D^b_c(\mcX_U,\Lambda)$ and assume $R^if_{U,*}F=0$. Then, for any $S\subsetneq B\setminus U$, the restriction morphism
\[
 H^i(\mcX_U,F)\longrightarrow\bigoplus_{b\in S}H^i(\xbhat^\circ,F)
\]
is surjective.
\end{lem}
\begin{proof}
Put $F_0=\tau^{\leq i}Rf_{U,*}F$. The assumption gives $F_0\in D^{\leq i-1}(U)$. Truncation and Leray identify the indicated restriction with
\[
 H^i(U,F_0)\longrightarrow
 \bigoplus_{b\in S}H^i(\wht B_b^\circ,F_0).
\]
Let $B_0=U\cup S\subsetneq B$, an affine open subcurve, and let $j:U\hookrightarrow B_0$. The comparison of Proposition~\ref{prop:MVsequence}, applied on the base to $j_!F_0$, gives
\begin{equation}\label{eq:base-MV}
\begin{split}
 H^i(U,F_0)\oplus\bigoplus_{b\in S}H^i(\wht B_b,j_!F_0)
 &\longrightarrow\bigoplus_{b\in S}H^i(\wht B_b^\circ,F_0)\\
 &\longrightarrow H^{i+1}(B_0,j_!F_0).
\end{split}
\end{equation}
The strictly henselian local schemes $\wht B_b$ have cohomological dimension zero, and the affine curve $B_0$ has prime-to-characteristic cohomological dimension at most one \cite[VI, Theorem 7.2]{Milne80}. Since $j_!F_0\in D^{\leq i-1}(B_0)$, both the local terms in degree $i$ and the final term in \eqref{eq:base-MV} vanish. This proves the surjectivity.
\end{proof}

\begin{cor}\label{cor:image}
Assume $R^{2d-1}f_{U,*}\Lambda=0$, and let $c\in B$ be any closed point with smooth fiber $\mcX_c$. After replacing $U$ by $U\setminus\{c\}$, the boundary in \eqref{eq:MVsequence} satisfies
\begin{equation}\label{eq:boundary-image}
 \operatorname{Im}(\partial)\subseteq
 \operatorname{Im}\left(H^{2d-2}(\mcX_c,\Lambda(d-1))
 \xrightarrow{i_{c,*}}H^{2d}(\mcX,\Lambda(d))\right).
\end{equation}
In particular, this image is algebraic if the curve classes span $H^{2d-2}(\mcX_c,\Lambda(d-1))$.
\end{cor}
\begin{proof}
Set $S=(B\setminus U)\setminus\{c\}$. Lemma~\ref{lem:surj} with $i=2d-1$ and $F=\Lambda(d)$ gives
\begin{equation}\label{eq:local-surjection}
 H^{2d-1}(\mcX_U,\Lambda(d))\twoheadrightarrow
 \bigoplus_{b\in S}H^{2d-1}(\xbhat^\circ,\Lambda(d)).
\end{equation}
The boundary in the sequence for $\mcX\setminus\mcX_c$ is consequently zero. By naturality of \eqref{eq:MVsequence}, the boundary image for $\mcX$ restricts to zero on this complement; the local contribution at $c$ disappears there. Localization and purity for the smooth divisor $i_c:\mcX_c\hookrightarrow\mcX$ give
\begin{equation}\label{eq:fiber-purity}
 H^{2d-2}(\mcX_c,\Lambda(d-1))\xrightarrow{i_{c,*}}
 H^{2d}(\mcX,\Lambda(d))\longrightarrow
 H^{2d}(\mcX\setminus\mcX_c,\Lambda(d)),
\end{equation}
which proves \eqref{eq:boundary-image}. The last assertion follows by compatibility with cycle pushforward.
\end{proof}

\subsection{\texorpdfstring{\textcolor{black}{Proof of the fibration theorem}}{Proof of the fibration theorem}}

\begin{proof}[Proof of Theorem~\ref{thm:WAKH2}]
The case $n=1$ is immediate, so assume $n>1$. Choose a nonempty open $U$ of smooth fibers with $c\notin U$, and put $B\setminus U=\{b_1,\ldots,b_s\}$. A global section $\sigma:B\to\mcX$ exists by \cite[Main Theorem]{deJongStarr_GHS}. The geometric generic fiber is simply connected, so finite-coefficient duality gives vanishing of its $H^{2d-1}(-,\Lambda(d))$. Smooth proper base change transports this vanishing over the connected smooth locus, hence
\begin{equation}\label{eq:odd-fiber-vanishing}
 R^{2d-1}\pi_{U,*}\Lambda(d)=0.
\end{equation}

Let $\alpha\in H^{2d}(\mcX,\Lambda(d))$, and write $\alpha_U$ and $\alpha_i$ for its restrictions to $\mcX_U$ and $\wht{\mcX}_{b_i}$. By Lemma~\ref{lem:local-classes}, each $\alpha_i$ is represented by $h_{i,*}[T_i]$ for a normal finite flat generically \'etale $R_i$-scheme of degree $D_i$. Equation~\eqref{eq:local-intersections} and the fiber divisor give
\begin{equation}\label{eq:degree-congruence}
 D_i\equiv\sum_\lambda m_{i\lambda}\theta_{b_i}(\alpha_i)_\lambda
 =\operatorname{Tr}_{\mcX}(\alpha\cup\pi^*[b_i])=:r\pmod n.
\end{equation}
All closed points of $B$ have the same cohomology class, so $r$ is independent of $i$.
Choose a positive $D\geq\max_iD_i$ with $D\equiv r\pmod n$. Adjoin to $T_i$ precisely $D-D_i$ disjoint copies of $\SP R_i$, mapping by $\sigma$. This is a multiple of $n$ copies, so it does not change the represented class. We now have normal finite flat generically \'etale sources of common positive degree $D$; the maps need not be embeddings.

By Lemma~\ref{lem:globalization}, there is a smooth integral projective cover $q:C\to B$ of degree $D$ with these completed local branches, as in \eqref{eq:globalized-covers}. Their maps to $\mcX$ are formal sections of $\mcX_C\to C$. Weak approximation after this finite base change supplies a section approximating all these maps to the orders required by Lemma~\ref{lem:map-continuity}. At a branch of ramification index $e$, congruence to order at least $eN$ in its uniformizer implies congruence modulo $t_i^N$. Properness extends the resulting rational section over all of $C$. Write $f:C\to\mcX$ for the composite. Lemma~\ref{lem:map-continuity} and proper pushforward give
\begin{equation}\label{eq:local-matching}
 \operatorname{cl}(f_*[C])|_{\wht{\mcX}_{b_i}}=\alpha_i
 \qquad(1\leq i\leq s).
\end{equation}

The finite-coefficient argument of \cite[Lemma 6.8]{Hasse0cycle} gives the trace isomorphism
\begin{equation}\label{eq:affine-trace}
 H^{2d}(\mcX_U,\Lambda(d))\xrightarrow{\sim}H^0(U,\Lambda)=\Lambda,
 \qquad\operatorname{cl}(\sigma|_U)\longmapsto1.
\end{equation}
Indeed, affine-curve cohomological dimension leaves only columns $0,1$ in Leray, \eqref{eq:odd-fiber-vanishing} eliminates the preceding row, and the top relative cohomology is $R^{2d}\pi_{U,*}\Lambda(d)=\Lambda$. Under \eqref{eq:affine-trace}, $\alpha_U$ maps to $r$, by its degree on any smooth fiber and \eqref{eq:degree-congruence}, while $\operatorname{cl}(f_*[C])|_{\mcX_U}$ maps to $D\equiv r$. 

\begingroup\color{black}
Together with \eqref{eq:local-matching}, the exact sequence \eqref{eq:MVsequence} therefore puts $\alpha-\operatorname{cl}(f_*[C])$ in its boundary image. By \eqref{eq:odd-fiber-vanishing} and Corollary~\ref{cor:image}, equation~\eqref{eq:boundary-image} supplies a class $\beta$ on the chosen smooth fiber with
\begin{equation}\label{eq:key}
 \alpha=\operatorname{cl}(f_*[C])+i_{c,*}\beta,
 \qquad\beta\in H^{2d-2}(\mcX_c,\Lambda(d-1)).
\end{equation}
 Equation~\eqref{eq:key} proves \eqref{eq:fibre-image}. Taking the cokernels of the cycle class maps in \eqref{eq:fibre-image} and tensoring with $\Lambda(1)$ gives, by \eqref{eq:kato-twist},
\begin{equation}\label{eq:fibre-defect}
 i_{c,*}:KH_2(\mcX_c,\Lambda)\twoheadrightarrow KH_2(\mcX,\Lambda).
\end{equation}
This is \eqref{eq:fibre-kato}; the maps agree because proper pushforward preserves cycle classes and commutes with the coefficient twist.
\endgroup

\end{proof}

\begingroup\color{black}
The conclusion for cycle classes follows from Theorem~\ref{thm:WAKH2}. We record its integral and $\ell$-adic forms.
\endgroup

\begingroup\color{black}
\begin{cor}\label{cor:WAITC}
    Let $\pi:\mcX\to B$ be a flat projective morphism of relative dimension $d\geq1$ over a smooth projective curve over an algebraically closed field $k$, with smooth separably rationally connected geometric generic fiber.
    Assume that $\mcX$ is smooth projective and that there is a smooth fiber $\mcX_c$ for which Conjecture~\ref{conj:IHC_RC} holds.
    Also assume that for every finite surjective morphism $C\to B$ from a smooth projective curve $C$, weak approximation holds for $\mcX_C=\mcX\times_B C\to C$.
    Then the cycle class map
    \[
    CH_1(\mcX)\otimes\ZZ_\ell\longrightarrow H^{2d}_{\text{\'et}}(\mcX,\ZZ_\ell(d))
    \]
    is surjective for every $\ell\ne\operatorname{char}k$.
    If $\mcX\to B$ is defined over the complex numbers, the cycle class map
    \[
    CH_1(\mcX)\longrightarrow H^{2d}_{\text{sing}}(\mcX,\ZZ)
    \]
    is surjective.
\end{cor}
\endgroup

\begingroup\color{black}
\begin{proof}
Smooth proper base change on the connected smooth locus transports vanishing of $H^{2d-1}(-,\ZZ/\ell^a)$ from the separably rationally connected geometric generic fiber to every smooth fiber, including $\mcX_c$. The inverse limit gives $H^{2d-1}(\mcX_c,\ZZ_\ell)=0$. The coefficient sequence consequently identifies
\begin{equation}\label{eq:fiber-coefficients}
 H^{2d-2}(\mcX_c,\ZZ_\ell(d-1))/\ell^a
 \simeq H^{2d-2}(\mcX_c,\ZZ/\ell^a(d-1)).
\end{equation}
The assumed $\ell$-adic cycle surjectivity on $\mcX_c$ and \eqref{eq:fiber-coefficients} imply its finite cycle surjectivity, in particular modulo $\ell$. Equation~\eqref{eq:fibre-image} of Theorem~\ref{thm:WAKH2} then gives finite cycle surjectivity on $\mcX$, since its fiber term is algebraic. Lemma~\ref{lem:reduction_mod_l}, with $\dim\mcX=d+1$, gives the asserted $\ell$-adic surjectivity; its finite-to-integral direction requires no additional odd-degree torsion-freeness assumption on $\mcX$.

Over $\CC$, comparison with singular cohomology and the same finite odd-degree cohomology vanishing show that the integral hypothesis on $\mcX_c$ implies the finite hypothesis for every $n$. Put
\[
 A=\operatorname{coker}\bigl(CH_1(\mcX)\longrightarrow H^{2d}_{\text{sing}}(\mcX,\ZZ)\bigr).
\]
The coefficient injection $H^{2d}_{\text{sing}}(\mcX,\ZZ)/n\hookrightarrow H^{2d}_{\text{sing}}(\mcX,\ZZ/n)$ and the finite surjectivity furnished by \eqref{eq:fibre-image} imply $A/nA=0$ for every $n$. Since $A$ is a finitely generated abelian group, it is zero. This proves the integral assertion, including torsion.
\end{proof}
\endgroup

\section{Quotients by finite groups}\label{sec:quotients}

\begingroup\color{black}
The following construction will play an important role in this section.

\begin{constr}\label{construction}
Let $X$ be a smooth projective variety over an algebraically closed field $k$, and let $G=\ZZ/\ell\ZZ$ be a cyclic group of order $\ell$, invertible in $k$, acting on $X$. Denote the quotient by $Z_0=X/G$. Choose a primitive $\ell$-th root of unity $\xi$, and let a generator $\sigma$ act on $\PP^1$ by
\[
 [X_0,X_1]\longmapsto[X_0,\xi X_1].
\]
Put $W=(X\times\PP^1)/G$ for the diagonal action, and choose a projective resolution $\tilde W\to W$ which is an isomorphism over the smooth locus, as in Remark~\ref{rem:G-resolution}.
\end{constr}
The following lemma summarizes some basic properties of this construction.
\begin{lem}\label{lem:constructionproperty}
The varieties $Z_0,W$ in Construction~\ref{construction} satisfy the following properties.
\begin{enumerate}
 \item $Z_0\times\PP^1$ is birational to $W$ over $Z_0$.
 \item The natural morphism $\tilde W\to W\to B=\PP^1/G$ exhibits $\tilde W,W$ as fibrations over $B$ whose general fibers are smooth projective varieties isomorphic to $X$.
\end{enumerate}
\end{lem}
\begin{proof}
Put $L=k(X)$, let $H\subset G$ be the kernel of the action on $L$, and set $h=|H|$, $Q=G/H$, and $F=L^Q$. For the affine coordinate $t$ with $\sigma(t)=\xi t$, one has $L(t)^H=L(t^h)$. Hilbert~90 for the cyclic Galois extension $L/F$ gives $a\in L^*$ satisfying $\sigma(a)=\xi^{-h}a$: the required norm is one, since $(\xi^{-h})^{\ell/h}=1$. With $z=at^h$, we obtain
\begin{equation}\label{eq:invariant-field}
 k(W)=L(t)^G=F(z)=k(X/G)(z).
\end{equation}
When $Q$ is trivial take $a=1$. Thus \eqref{eq:invariant-field} also proves the assertion for nonfaithful actions.
Over $B\setminus\{0,\infty\}$, the map $\PP^1\to B$ is a $G$-torsor. Descent of the smooth product family gives a smooth isotrivial family with geometric fibers $X$. The resolution in Construction~\ref{construction} preserves this smooth locus.
\end{proof}
\endgroup

\begingroup\color{black}
\begin{lem}\label{lem:tame-isotrivial-WA}
Suppose $X$ is separably rationally connected and $G$ has prime order $\ell\ne\operatorname{char}k$. For every finite surjective morphism $C\to B$ of smooth projective curves, including inseparable morphisms, the generic fiber of $\tilde W\times_B C\to C$ satisfies weak approximation.
\end{lem}
\endgroup

\begin{proof}
\begingroup\color{black}
Write $s=t^\ell$ on $B$ and put $K=k(C)$. The splitting algebra
\begin{equation}\label{eq:kummer-basechange}
 K[t]/(t^\ell-s)
\end{equation}
is \'etale. Since $\ell$ is prime and $K$ contains the $\ell$-th roots of unity, it is split or a cyclic extension of degree $\ell$, even when $C\to B$ is inseparable. Normalize $C$ in a field factor and denote the resulting smooth projective curve by $C'$. In the nonsplit case this is a tame cyclic cover. The generic fiber becomes the constant variety $X$ over $k(C')$, with smooth proper model $X\times C'$. At a ramified point the completed cover has parameters
\begin{equation}\label{eq:local-inertia}
 k[[u]]\subset k[[v]],\qquad u=v^\ell,\qquad
 \sigma(v)=\xi_c v,
\end{equation}
for a suitable primitive $\ell$-th root $\xi_c$. Thus the local model upstairs is a smooth separably rationally connected product and the inertia is tame cyclic. At unramified places it has good reduction.
\endgroup


We use the argument of \cite[\S\S2--4]{WAIsotrivial}. 
We explain the tame extension of that argument and in particular address an issue in characteristic $2$ overlooked in the paper.  The key theorem is the existence of an equivariant rational curve connecting any pair of fixed points (Theorem~3.1). The proof works in any characteristic other than $2$ as long as the order of the cyclic group is invertible over the field.

In characteristic two, the chosen pencil in Proposition 3.6 of \cite{WAIsotrivial} contains a member $X_0^{2\ell}-X_1^{2\ell}$ that is nonreduced. 
We just need to fix this.
Choose distinct $\lambda,\mu\in k^*$. On $\PP^1_X\times\PP^1_Y$, use the pencil generated by
\begin{align}
 C_0&:(X_1^\ell-\lambda X_0^\ell)(X_1^\ell-\mu X_0^\ell)Y_0Y_1=0,
 \label{eq:pencil-zero}\\
 C_\infty&:(X_0^\ell Y_1+X_1^{\ell} Y_0)
 \bigl(X_0^\ell Y_1+((\lambda+\mu)X_1^\ell-\lambda\mu X_0^\ell)Y_0\bigr)=0.
 \label{eq:pencil-infty}
\end{align}
The action scales $X_1$ by $\xi$ and fixes $Y_0,Y_1$. The two components of \eqref{eq:pencil-infty} are smooth rational. 
Their intersection equation is $(X_1^\ell -\lambda X_0^\ell)(X_1^\ell-\mu X_0^\ell)=0$, giving two free orbits of $\ell$ distinct points, where they meet transversely. Its first component contains the fixed markings
\[
 ([1,0],[1,0]),\qquad ([0,1],[0,1]).
\]
The divisor \eqref{eq:pencil-zero} is nodal, with two horizontal components and two free orbits of vertical components, and all its nodes have trivial stabilizers. The section defining \eqref{eq:pencil-infty} is nonzero at every node of \eqref{eq:pencil-zero}; therefore the pencil smooths those nodes. Away from the nodes smoothness follows by openness and properness. At the two markings, the linear coefficients of a member $aC_0+bC_\infty$ are nonzero on $a-b\ne0$ and $a+b\ne0$, respectively. Both pencil endpoints belong to this open, also in characteristic two. The tangent characters there are $\xi$ and $\xi^{-1}$. This supplies exactly the required degeneration with two distinct free orbits, replacing the printed orbits $t^\ell=1,-1$ by $t^\ell=\lambda,\mu$. The remaining equivariant deformation and jet-induction arguments are unchanged.
Thus the tame argument applies to every cover \eqref{eq:kummer-basechange}, proving the lemma.
\end{proof}

\begingroup\color{black}
\begin{prop}\label{prop:cyclic}
Let a cyclic group $H$ of prime order $r\ne\operatorname{char}k$ act on a smooth projective separably rationally connected variety $X$. Let $Y$ be a smooth projective model of $X/H$. For every $n$ invertible in $k$, put $\Lambda=\ZZ/n$. The rational quotient map induces a surjection
\[
 q_*:KH_2(X,\Lambda)\twoheadrightarrow KH_2(Y,\Lambda).
\]
\end{prop}
\begin{proof}
Dimension zero is immediate. If the action is trivial, $q$ is birational and Lemma~\ref{lem:graphs} applies. Otherwise apply Construction~\ref{construction} with $G=H$ and $\ell=r$, and write $\mcT=\tilde W$. The morphism $\mcT\to B$ is flat, since $\mcT$ is integral and $B$ is a smooth curve. By Lemma~\ref{lem:constructionproperty}, its fibers over $B^\circ=B\setminus\{0,\infty\}$ are smooth and isomorphic to $X$. Lemma~\ref{lem:tame-isotrivial-WA} supplies weak approximation after every finite curve base change. Therefore Theorem~\ref{thm:WAKH2}, in the form \eqref{eq:fibre-defect}, gives
\begin{equation}\label{eq:cyclic-fibre}
 i_{c,*}:KH_2(X,\Lambda)\twoheadrightarrow KH_2(\mcT,\Lambda)
 \qquad(c\in B^\circ(k)).
\end{equation}
To specify this fiber identification, write $s=t^r$ on $B$ and choose $a\in k^*$ with $a^r=c$. The map $x\mapsto[(x,a)]$ identifies $X$ with $\mcT_c=W_c$. Under this identification, the projection $W\to X/H$ restricts to the quotient map $X\to X/H$ itself.

The invariant-field identity \eqref{eq:invariant-field} is over $k(X/H)$. It follows that the rational projection $\rho:\mcT\dashrightarrow Y$ corresponds, under a birational map $\mcT\dashrightarrow Y\times\PP^1$, to the first projection. The birational and projective-bundle assertions of Lemma~\ref{lem:graphs} give an isomorphism
\begin{equation}\label{eq:rho}
 \rho_*:KH_2(\mcT,\Lambda)\xrightarrow{\sim}KH_2(Y,\Lambda).
\end{equation}
The map $\rho$ is defined at the generic point of $\mcT_c$: it is a morphism on the inverse image of a dense open where $Y$ agrees with $X/H$. Its restriction is $q$. The restriction identity \eqref{eq:restriction} hence gives
\begin{equation}\label{eq:coefficient-one}
 \rho_*i_{c,*}=q_*:KH_2(X,\Lambda)\longrightarrow KH_2(Y,\Lambda).
\end{equation}
Combining the surjection \eqref{eq:cyclic-fibre}, the isomorphism \eqref{eq:rho}, and \eqref{eq:coefficient-one} proves the proposition. There is no factor $r$ in \eqref{eq:coefficient-one}: the chosen $X\times\{a\}$ maps isomorphically to the single quotient fiber $W_c$, and we do not push forward the union of all $r$ lifted fibers. Thus no division by the group order is involved.
\end{proof}
\endgroup

\begingroup\color{black}
\begin{prop}\label{prop:finite}
For $X,G,Y$ as in Theorem~\ref{thm:main} and $\Lambda=\ZZ/n$ with $n$ invertible in $k$, the rational quotient pushforward $q_*:KH_2(X,\Lambda)\to KH_2(Y,\Lambda)$ is surjective.
\end{prop}
\begin{proof}
Reduce to $n=\ell^a$ by Chinese Remainder Theorem. First let $P:=G$ be an
$\ell$-group, and induct on its order. The trivial group case is the birational assertion of Lemma~\ref{lem:graphs}. Choose a central subgroup
$H\subset P$ of order $\ell$.
Remark~\ref{rem:G-resolution} and \eqref{eq:equivariant-resolution} give a smooth projective separably rationally connected model $Z$ of $X/H$ carrying the residual $P/H$-action.
Proposition~\ref{prop:cyclic} and induction give surjections
\[
 KH_2(X,\Lambda)\xrightarrow{q_{H,*}}KH_2(Z,\Lambda)
 \xrightarrow{q_{P/H,*}}KH_2(Y_P,\Lambda)
\]
for a smooth projective tame quotient model $Y_P$ of $X/P$.
The invariant-field identity
\begin{equation}\label{eq:central}
 (k(X)^H)^{P/H}=k(X)^P
\end{equation}
and functoriality in Lemma~\ref{lem:graphs} identify the composite
with $q_{P,*}$. Birational invariance in Lemma~\ref{lem:graphs} transfers the assertion to every smooth projective model of $X/P$. This also covers subgroups acting trivially.

For general $G$, let $P$ be a Sylow $\ell$-group and $N$ the kernel subgroup of
the group action on $k(X)$. 
The preceding case gives surjectivity to
a smooth projective model $Z$ of $X/P$. The map
$f:Z\dashrightarrow Y$ has degree
\begin{equation}\label{eq:sylow}
 m=[k(X)^P:k(X)^G]=[G:PN]\mid[G:P],
 \qquad(m,\ell)=1.
\end{equation}
By \eqref{eq:pushpull}, $f_*f^*=m\operatorname{id}$, so $f_*$
is surjective. Functoriality gives
$q_{G,*}=f_*q_{P,*}$, proving the proposition.
We emphasize that we assume the smooth projective model $Y$ exists, and that only the $\ell$-group intermediate quotients require tame resolutions (See Remark \ref{rem:G-resolution}).
\end{proof}
\endgroup

\begingroup\color{black}
\begin{proof}[Proof of Theorem~\ref{thm:main}]
For every $g\in G$, the rational maps $q\circ g$ and $q$ agree.
Lemma~\ref{lem:graphs} gives
\begin{equation}\label{eq:coinvariants}
 q_*g_*=q_*.
\end{equation}
Proposition~\ref{prop:finite} therefore factors as a surjection
$KH_2(X,\Lambda)_G\twoheadrightarrow KH_2(Y,\Lambda)$.
This proves \eqref{eq:main}, with
the stated birational-model compatibility supplied by
Lemma~\ref{lem:graphs}.

For the integral assertion put
\begin{equation}\label{eq:integral-defect}
 KH_2(V,\ZZ)=\operatorname{coker}\bigl(CH_1(V)\to H_2(V,\ZZ)\bigr),
\end{equation}
using \cite[equation (3.1)]{zerocycleLaurent}.
The integral versions of Lemmas~\ref{lem:boundary} and~\ref{lem:graphs}
hold over $\CC$, using a degree-one resolution of the proper source
support instead of a prime-to-$\ell$ alteration. Thus the integral
quotient graph map is defined and satisfies \eqref{eq:coinvariants}.
By Remark~\ref{rem:G-resolution}, both $X$ and $Y$ are rationally connected, hence have $H_1(-,\ZZ)=0$.
Universal coefficients and compatibility of cycle classes imply
\begin{equation}\label{eq:integral-reduction}
 H_2(V,\ZZ)/n\simeq H_2(V,\ZZ/n),\qquad
 KH_2(V,\ZZ)/n\simeq KH_2(V,\ZZ/n)\quad(V=X,Y),
\end{equation}
where finite-coefficient singular homology is identified with \'etale
Borel--Moore homology by singular--\'etale comparison; compare
\cite[Lemma 3.3]{zerocycleLaurent}.
Let $C$ be the cokernel of the integral pushforward
$KH_2(X,\ZZ)\to KH_2(Y,\ZZ)$.
Proposition~\ref{prop:finite} and \eqref{eq:integral-reduction}
give $C/nC=0$ for every $n>0$. It is a finitely generated abelian
group, being a quotient of singular $H_2(Y,\ZZ)$. Therefore $C=0$.
Finally \eqref{eq:coinvariants} and \eqref{eq:integral-defect}
give \eqref{eq:integral-main}.
\end{proof}
\endgroup

\begingroup\color{black}
\begin{proof}[Proof of Corollary~\ref{cor:quotientITC}]
By Remark~\ref{rem:G-resolution}, the supplied smooth projective model $Y$ is separably rationally connected. For a smooth proper separably rationally connected $d$-fold, simple connectivity and finite-coefficient duality give
\begin{equation}\label{eq:SRC-odd-cohomology}
 H^{2d-1}(-,\ZZ/\ell^a)=0
 \qquad(\ell\ne\operatorname{char}k,\ a\geq1).
\end{equation}
The assumed $\ell$-adic cycle surjectivity on $X$, \eqref{eq:SRC-odd-cohomology}, and the coefficient sequence as in \eqref{eq:fiber-coefficients} give $KH_2(X,\ZZ/\ell^a)=0$. Equations~\eqref{eq:kato-twist} and \eqref{eq:main} of Theorem~\ref{thm:main} then give $KH_2(Y,\ZZ/\ell^a)=0$. By \eqref{eq:defect}, this is finite cycle surjectivity on $Y$. Lemma~\ref{lem:reduction_mod_l} proves the $\ell$-adic assertion.

Over $\CC$, the hypothesis says $KH_2(X,\ZZ)=0$ in \eqref{eq:integral-defect}. The integral surjection \eqref{eq:integral-main} of Theorem~\ref{thm:main} gives $KH_2(Y,\ZZ)=0$, which is the required integral cycle surjectivity, including torsion. Finally, for $X=\PP^d$ the source cohomology is algebraic by the projective bundle formula. If $|G|$ is invertible in $k$, \eqref{eq:equivariant-resolution} supplies a smooth projective quotient model, proving the last assertion.
\end{proof}

\endgroup

\bibliographystyle{alpha}
\bibliography{MyBib}

\end{document}